\documentclass{amsart}
\usepackage{enumerate,color,amsmath,amssymb,amsfonts,bm}
\usepackage{frcursive,fullpage,graphicx,psfrag}
\usepackage{color}

\newtheorem{Theorem}{Theorem}[section]
\def\bbone{{\mathchoice {\rm 1\mskip-4mu l} {\rm 1\mskip-4mu l}
{\rm 1\mskip-4.5mu l} {\rm 1\mskip-5mu l}}}

\begin{document}

\author{Abdelmalek Abdesselam}
\address{Abdelmalek Abdesselam,
Department of Mathematics,
P. O. Box 400137,
University of Virginia,
Charlottesville, VA 22904-4137, USA}
\thanks{A.A. was supported in part by the National Science
Foundation under grant DMS \# 0907198.}
\email{malek@virginia.edu}

\author{Ajay Chandra}
\address{Ajay Chandra,
Department of Mathematics,
P. O. Box 400137,
University of Virginia,
Charlottesville, VA 22904-4137, USA}
\email{ac2yx@virginia.edu}

\author{Gianluca Guadagni}
\address{Gianluca Guadagni,
Department of Mathematics,
P. O. Box 400137,
University of Virginia,
Charlottesville, VA 22904-4137, USA}
\email{gg5d@virginia.edu}

\title{Rigorous quantum field theory functional integrals over the $p$-adics: research announcement}

\begin{abstract}
In this short note we announce the construction of scale invariant non-Gaussian generalized stochastic
processes over three dimensional $p$-adic space. The construction includes that of the associated squared field and
our result shows this squared field has a dynamically generated anomalous dimension which rigorously
confirms a prediction made more than forty years ago, in an essentially identical situation, by Kenneth G. Wilson. 
We also prove a mild form of universality for the model under consideration.
\end{abstract}

\maketitle

\section{Preamble}
The findings presented in this brief announcement can be construed as a direct continuation
of a line of investigation initiated and developed by the school of mathematical thought
around Andrei N. Kolmogorov, with particular regard to scale invariant generalized stochastic
processes~\cite{Kolmogorov,Minlos,GelfandS,GelfandV,Sinai,Dobrushin1,Dobrushin2}.
Such processes typically arise as scaling limits of models of statistical mechanics (see, e.g., the recent work~\cite{CamiaGN}).
For simplicity and as a stepping stone towards the study of such processes over real space, we consider here a natural analogue
over three dimensional $p$-adic space where harmonic analysis is more forgiving.
In the next two sections we will introduce the definitions needed for the formulation of our theorem which is given in the fourth and last section.
For a more detailed presentation of real or complex valued $p$-adic analysis which is suitable for a probability theory or mathematical physics audience, the reader
may consult~\cite{AlbeverioKS} and~\cite{VladimirovVZ}. An additional suitable and recommended reference is~\cite{GoldfeldH} which makes the pedagogically excellent choice of focusing on completions of the field $\mathbb{Q}$ and thus renders the material accessible without extensive knowledge of algebraic number theory.

\section{Generalities about $p$-adics}
Let $p$ be a prime number and consider the $p$-adic absolute value $|\cdot|_p$ on $\mathbb{Q}$
defined by $|x|_p=0$ if $x=0$ and $|x|_p=p^{-k}$ if $x=\frac{a}{b}\times p^k$ where $a,k\in\mathbb{Z}$ and $b$, a positive integer, 
are such that $a,b$ are coprime and neither are divisible by $p$.
The field $\mathbb{Q}_p$ of $p$-adic numbers is the completion of $\mathbb{Q}$ with respect to this absolute value.
Every $x$ in $\mathbb{Q}_p$
has a unique convergent series representation
\[
x=\sum_{n\in\mathbb{Z}} a_n p^n
\]
where the digits $a_n$ belong to $\{0,1,\ldots,p-1\}$ and at most finitely many of them are nonzero for negative $n$.
The absolute value of $x\neq 0$ can be recovered from this representation as $|x|_p=p^{-v_p(x)}$ where
\[
v_p(x)=\min\{n\in\mathbb{Z}\ |\ a_n\neq 0\}\ .
\]
Using the same representation one can define the fractional (or polar) part of $x$ which is
$\{x\}_p=\sum_{n<0} a_n p^n$.
The closed unit ball $\mathbb{Z}_p=\{x\in\mathbb{Q}_p|\ |x|\le 1\}$ is a compact subring of $\mathbb{Q}_p$.
From now on we will drop the $p$ subscript from the absolute value.
The additive Haar measure on $\mathbb{Q}_p$ normalized so that $\mathbb{Z}_p$ has measure one
will simply be denoted by ${\rm d}x$.
In $d$ dimensions,
the $p$-adic norm of a vector $x=(x_1,\ldots,x_d)\in\mathbb{Q}_p^d$ is defined as
$|x|=\max\{|x_1|,\ldots,|x_d|\}$. The product measure ${\rm d}^d x$ obtained from the previous one-dimensional measure
is invariant by the subgroup $GL_d(\mathbb{Z}_p)$ of $GL_d(\mathbb{Q}_p)$.
The subgroup $GL_d(\mathbb{Z}_p)$ is defined as the set of $d\times d$ matrices which together with their inverses have entries in $\mathbb{Z}_p$. This subgroup is the maximal compact subgroup of $GL_d(\mathbb{Q}_p)$ (unique up to conjugacy) and is the
natural analogue of the orthogonal group $O(d)$ acting on $\mathbb{R}^d$.
The use of the maximum in the definition of the norm is motivated by the resulting invariance with respect to $GL_d(\mathbb{Z}_p)$.

The space of real (resp. complex) test functions $S(\mathbb{Q}_p^d,\mathbb{R})$ (resp. $S(\mathbb{Q}_p^d,\mathbb{C})$)
is the Schwartz-Bruhat space of compactly supported locally constant real-valued (resp. complex-valued) functions on $\mathbb{Q}_p^d$.
If we do not specify the target, then we mean $\mathbb{R}$.
Recall that a seminorm on $S(\mathbb{Q}_p^d)$ is a function $N:S(\mathbb{Q}_p^d)\rightarrow [0,\infty)$
which satisfies the usual norm axioms except the requirement that $N(f)=0$ implies $f=0$.
The coarsest topology on $S(\mathbb{Q}_p^d)$ which makes all possible seminorms continuous is called the finest locally convex
topology and it is the one we use. The space of distributions $S'(\mathbb{Q}_p^d)$
simply is the topological dual of $S(\mathbb{Q}_p^d)$ which turns out to be the algebraic dual.
The probability measures we will be interested in are defined on the measurable space $(S'(\mathbb{Q}_p^d),\mathcal{C})$
where $\mathcal{C}$ is the cylinder $\sigma$-algebra, i.e., the smallest which makes all coordinate functions $\phi\mapsto \phi(f)$,
indexed by $f\in S(\mathbb{Q}_p^d)$, measurable.
The Fourier transform of a complex valued test function $f$ is defined by
\[
\widehat{f}(k)=\int_{\mathbb{Q}_p^d} f(x)\exp(-2i\pi\{k\cdot x\}_p)\ {\rm d}^d x
\]
where $k\cdot x=k_1 x_1+\cdots+k_d x_d$ and the rational $\{k\cdot x\}_p$ is seen as a real number.
One has that the characteritic function of $\mathbb{Z}_p^d$ is fixed by the Fourier transform, that is
$\widehat{\bbone}_{\mathbb{Z}_p^d}=\bbone_{\mathbb{Z}_p^d}$.
From this it easily follows that the space $S(\mathbb{Q}_p^d,\mathbb{C})$ is stable by Fourier transform.
One can also define the Fourier transform of distributions by duality.
The analogue of the nuclear theorem holds in this setting as well as that of the Bochner-Minlos Theorem~\cite{Minlos,GelfandV}.

If one views a point $x$ in $\mathbb{Q}_p^d$ as a column vector
then one has a left-action of $GL_d(\mathbb{Z}_p)$ on points simply by matrix multiplication.
It results in left-actions on test functions $f$, distributions 
$\phi$ and more generally $n$-linear forms $W$ (which are automatically continuous)
on $S(\mathbb{Q}_p^d)$, using
\[
(M\cdot f)(x)=f(M^{-1}x)\ ,
\]
\[
(M\cdot\phi)(f)=\phi(M^{-1}\cdot f)\ ,
\]
\[
(M\cdot W)(f_1,\ldots,f_n)=W(M^{-1}\cdot f_1,\ldots,M^{-1}\cdot f_n)\ .
\]
Such objects are called rotation invariant if they are preserved by all $M\in GL_d(\mathbb{Z}_p)$.
If one formally thinks of a distribution $\phi$ as a `function' via
\[
\phi(f)=\int_{\mathbb{Q}_p^d} \phi(x) f(x)\ {\rm d}^d x
\]
then the choice of definition means ``$(M\cdot\phi)(x)=\phi(M^{-1}x)$''.
Thus, a distribution $\phi$ is rotation invariant if ``$\phi(M^{-1}x)=\phi(x)$'' for all $M$ and $x$. 
A probability measure $\nu$ is called rotation invariant if for any $M\in GL_d(\mathbb{Z}_p)$
the push-forward (or direct image) of $\nu$ by the map $\phi\mapsto M\cdot\phi$ is $\nu$ itself.
This means the last formal equality between quotes holds in distribution: both sides have the same probability law.

Likewise regarding translations, one can define for $y\in\mathbb{Q}_p^d$ the transformations
\[
\tau_y(x)=x+y\ ,
\]
\[
\tau_y(f)(x)=f(x-y)\ ,
\]
\[
\tau_y(\phi)(f)=\phi(\tau_{-y}(f))\ ,
\]
\[
\tau_y(W)(f_1,\ldots,f_n)=W(\tau_{-y}(f_1),\ldots,\tau_{-y}(f_n))\ .
\]
One then defines the notion of invariance by translation for $n$-linear forms and probability measures on $S'(\mathbb{Q}_p^d)$
in the same way as before.

We now turn to scaling transformations. Given $\lambda\in\mathbb{Q}_p^\ast=\mathbb{Q}_p\backslash\{0\}$, we write:
\[
(\lambda\cdot f)(x)=f(\lambda^{-1}x)\ ,
\]
\[
(\lambda\cdot\phi)(f)=|\lambda|^d\ \phi(\lambda^{-1}\cdot f)\ .
\]
This corresponds to the formal equation ``$(\lambda\cdot \phi)(x)=\phi(\lambda^{-1}x)$''.
A distribution $\phi$ is called partially scale invariant with homogeneity $\alpha\in\mathbb{R}$ with respect to a subgroup
$H$ of the full scaling group $p^{\mathbb{Z}}\subset \mathbb{Q}_p^\ast$ if $\lambda\cdot\phi=|\lambda|^{-\alpha}\phi$
for all $\lambda\in H$.
This formally means ``$|\lambda|^{\alpha}\phi(\lambda^{-1} x)=\phi(x)$''.
A probability measure on $S'(\mathbb{Q}_p^d)$ is called partially scale invariant with homogeneity $\alpha$
with respect to the subgroup $H$ if the push-forward of $\nu$ by
the map $\phi\mapsto |\lambda|^{\alpha}(\lambda\cdot\phi)$
is $\nu$ itself, for all $\lambda\in H$.
Namely, the last equation between quotes is to be interpreted as saying that the random fields on both sides have the same law.

\section{Specifics about the model under consideration}
Now let us pick $d=3$ and for $0<\epsilon<1$ let us denote the quantity $\frac{3-\epsilon}{4}$ by the symbol $[\phi]$.
Let $L=p^l$ for some integer $l\ge 1$. For $r\in\mathbb{Z}$ (typically negative), we consider
the bilinear form on $S(\mathbb{Q}_p^3)$ given by
\[
C_r(f,g)=\int_{\mathbb{Q}_p^3} 
\frac{\widehat{f}(-k)\widehat{g}(k)\bbone\{|k|\le L^{-r}\}}{|k|^{3-2[\phi]}}\ {\rm d}^3 k
\]
where we use $\bbone\{\cdots\}$ for the characteristic function of the condition between braces.
By the $p$-adic version of the Bochner-Minlos Theorem, there is a unique probability measure $\mu_{C_r}$
on $S'(\mathbb{Q}_p^3)$ such that for any $f\in S(\mathbb{Q}_p^3)$
\[
\left\langle e^{i\phi(f)}\right\rangle_{\mu_{C_r}}=\exp\left(-\frac{1}{2}C_r(f,f)\right)
\]
where we used the statistical mechanics notation for the expectation with respect to $\phi$
sampled according to the measure $\mu_{C_r}$.
Note that one can write, with a slight abuse of notation
\[
C_r(f,g)=\int_{\mathbb{Q}_p^{3\times 2}}
\ C_r(x-y)f(x)g(y)\ {\rm d}^3 x\ {\rm d}^3 y
\]
where the function $C_r$ is explicitly given by
\[
C_r(x)=\sum_{n=lr}^{\infty}
p^{-2n[\phi]}\left[
\bbone_{\mathbb{Z}_p^3}(p^n x)-p^{-3}\bbone_{\mathbb{Z}_p^3}(p^{n+1} x)
\right]\ .
\]

The measure $\mu_{C_r}$ is supported on distributions given by bonafide functions which are locally constant at scale $L^r$, namely,
constant on each coset in $\mathbb{Q}_p^3/(L^{-r}\mathbb{Z}_p)^3$, the latter quotient playing the role of the lattice of mesh $L^r$.
Note that since $|p|=p^{-1}$ where the $p$ on the left wears its $p$-adic hat while the one on the right is viewed as a real number,
the volume of $(L^{-r}\mathbb{Z}_p)^3$ is $L^{3r}$ in accordance with the intuitive image of a three-dimensional
box with linear dimension $L^r$.
For $s\in\mathbb{Z}$ (typically positive), we use the notation $\Lambda_s=\{x\in\mathbb{Q}_p^3|\ |x|\le L^s\}$ and we also define the Wick
powers
\[
:\phi^2:_{C_r}(x)=\phi(x)^2-C_r(0)\ ,
\]
\[
:\phi^4:_{C_r}(x)=\phi(x)^4-6\ C_r(0)\ \phi(x)^2+3\ C_r(0)^2
\]
and, given $g>0$ as well as $\mu\in\mathbb{R}$, the potential
\[
V_{r,s}(\phi)=\int_{\Lambda_s}
\left\{
L^{-(3-4[\phi])r}\ g :\phi^4:_{C_r}(x)+L^{-(3-2[\phi])r}\ \mu:\phi^2:_{C_r}(x)
\right\}\ {\rm d}^3 x\ .
\]
By the previous remarks, the measure 
\[
{\rm d}\nu_{r,s}(\phi)=\frac{1}{\mathcal{Z}_{r,s}}e^{-V_{r,s}(\phi)}{\rm d}\mu_{C_r}(\phi)
\]
is a well defined probability measure on $S'(\mathbb{Q}_p^3)$ with finite moments.
The normalization factor $\mathcal{Z}_{r,s}$ is at least equal to one as can be seen from Jensen's inequality.
We will denote expectations with respect to $\nu_{r,s}$ by $\langle\cdots\rangle_{r,s}$.
Finally, given a locally constant $\phi$ at scale $L^r$,
we define an element $N_{r}[\phi^2]$ of $S'(\mathbb{Q}_p^3)$
by letting it act on $j\in S(\mathbb{Q}_p^3)$ via
\[
N_{r}[\phi^2](j)=Z_2^{\ r} \int_{\mathbb{Q}_p^3}
\left(
Y_2 :\phi^2:_{C_r}(x)-Y_0\ L^{-2r[\phi]}
\right)\ j(x)\ {\rm d}^3 x
\]
where $Z_2,Y_0, Y_2$ are parameters used in the construction.

We will also use the notation
\[
\bar{g}_\ast=\frac{(p^\epsilon-1)}{36\ L^{\epsilon}(1-p^{-3})}\ .
\]

\section{The result}

\begin{Theorem}
\label{mainthm}
\ 

$\exists \rho>0$, $\exists L_0$, $\forall L\ge L_0$,
$\exists \epsilon_0>0$, $\forall\epsilon\in (0,\epsilon_0]$,
one can find $\eta_{\phi^2}>0$ and functions $\mu(g),Y_0(g),Y_2(g)$
of $g$ in the interval $(\bar{g}_\ast-\rho\epsilon^{\frac{3}{2}},\bar{g}_\ast+\rho\epsilon^{\frac{3}{2}})$,
such that if one sets $\mu=\mu(g)$, $Z_2=L^{-\frac{1}{2}\eta_{\phi^2}}$, $Y_0=Y_0(g)$
and $Y_2=Y_2(g)$ in the previous definitions, then for all collections of test functions
$f_1,\ldots,f_n,j_1,\ldots,j_m$, the limits
\[
\lim_{\substack{r\rightarrow -\infty\\ s\rightarrow\infty}}
\left\langle
\phi(f_1)\cdots\phi(f_n)
N_r[\phi^2](j_1)\cdots N_r[\phi^2](j_m)
\right\rangle_{r,s}
\] 
exist and do not depend on the
order in which the $r\rightarrow -\infty$ and $s\rightarrow\infty$ limits are taken.
Moreover, the resulting quantities or correlators henceforth similarly and formally denoted by dropping the $r$ and $s$ subscripts
(and using squares, 4-th powers, etc., for repeats)
satisfy the following properties:

1) They are invariant by translation and rotation.

2) They satisfy the partial scale invariance property
\[
\left\langle
\phi(\lambda\cdot f_1)\cdots\phi(\lambda\cdot f_n)
\ N[\phi^2](\lambda\cdot j_1)\cdots N[\phi^2](\lambda\cdot j_m)
\right\rangle=\qquad\qquad\qquad\qquad
\]
\[
\qquad\qquad\qquad|\lambda|^{(3-[\phi])n+(3-2[\phi]-\frac{1}{2}\eta_{\phi^2})m}
\left\langle
\phi(f_1)\cdots\phi(f_n)
\ N[\phi^2](j_1)\cdots N[\phi^2](j_m)
\right\rangle
\]
for all $\lambda\in L^{\mathbb{Z}}$.

3) They satisfy the nontriviality conditions
\[
\langle \phi(\bbone_{\mathbb{Z}_p^3})^4
\rangle - 3\langle \phi(\bbone_{\mathbb{Z}_p^3})^2
\rangle <0\ ,
\]
\[
\langle N[\phi^2](\bbone_{\mathbb{Z}_p^3})^2
\rangle=1\ .
\]

4) The pure $\phi$ correlators are the moments of a unique probability measure $\nu_{\phi}$ on $S'(\mathbb{Q}_p^2)$
with finite moments. This measure
is translation and rotation invariant. It is also partially scale invariant with homogeneity $-[\phi]$ with respect to
the scaling subgroup $L^{\mathbb{Z}}$.

5) The pure $N[\phi^2]$ correlators are the moments of a unique probability measure $\nu_{\phi^2}$ on $S'(\mathbb{Q}_p^2)$
with finite moments. This measure
is translation and rotation invariant. It is also partially scale invariant with homogeneity $-2[\phi]-\frac{1}{2}\eta_{\phi^2}$
with respect to the scaling subgroup $L^{\mathbb{Z}}$.

6) The measures $\nu_\phi$ and $\nu_{\phi^2}$ satisfy a mild form of universality: they do not depend on $g$ in the above-mentioned interval.
\end{Theorem}

A few comments are in order.
\begin{itemize}
\item
The proof of our theorem is quite long and details will appear in a forthcoming paper~\cite{ACG1}.
We believe this proof may serve as a blueprint
that could guide future efforts devoted to the proof of a similar result over the reals.

\item
A trivial consequence of the sign symmetry $\phi\rightarrow -\phi$ is that $\langle\phi(f)\rangle$ is identically zero.
Almost as trivial, in view of the mentioned invariance properties, is that also $\langle N[\phi^2](j)\rangle$ is identically zero.
In other words $N[\phi^2]$ shares similar properties with the Wick normal product $:\phi^2:$ in the free, i.e., Gaussian case.

\item
We emphasize that $Z_2$ is independent of
$g$. However the quantities $\mu,Y_0,Y_2$ are allowed to depend on $g$. This is to be born in mind for the proper
interpretation of the universality property.
However, it probably did not escape the attentive reader that our theorem leaves the door open to the dependence
of the anomalous dimention $\eta_{\phi^2}$ on $L$. In fact our proof also produces the small epsilon asymptotic
$\eta_{\phi^2}=\frac{2}{3}\epsilon+o(\epsilon)$ which shows the independence on $L$ at least to first order.
We are making rapid progress towards proving that $\eta_{\phi^2}$ is indeed independent of $L$
as well as removing the restriction to the subgroup
$L^{\mathbb{Z}}$ regarding the mentioned scale invariance properties. We do not yet claim this stronger result
as a theorem since some details of our argument~\cite{ACG2} still remain to be checked.

\item
The absence of anomalous scaling for the field $\phi$ together with the presence of such anomalous scaling for the $\phi^2$ field
is a prediction made by K. G. Wilson for a very similar model in~\cite{Wilson}. Therefore our theorem provides a rigorous
mathematical proof of this statement.

\end{itemize}

\end{document}